\input amstex
\documentstyle{amsppt}
\input epsf
\magnification 1200
\NoBlackBoxes
\vcorrection{-1cm}

\def\letterwise{\equiv}

\def\hA{\hat{\Cal A}}
\def\wt{\operatorname{wt}}
\def\simh{\sim_{\text{h}}}

\rightheadtext{Criterion for Hurwitz equivalence for 3-braid}
\topmatter
\title
      Criterion of Hurwitz equivalence for quasipositive factorizations of 3-braids
\endtitle
\author
            S.~Yu.~Orevkov
\endauthor
\address   Steklov Math. Inst. (Moscow) and
           Univ. Paul Sabati\'e (Toulouse-3)
\endaddress
\email orevkov\@math.ups-tlse.fr
\endemail
\endtopmatter

\vskip -5mm

\document
The problem of Hurwitz equivalence of $n$-tuples of braids
appears in the study of braid monodromy of algebraic curves in $\Bbb C^2$.
It was considered by many authors, see [1, 3, 4]
and references in [3]. We give an answer for the case mentioned in the title.

Let 
$\bold B_3=\langle\Cal A\mid
\sigma_2\sigma_1=\sigma_1\sigma_0=\sigma_0\sigma_2\rangle$,
$\Cal A=\{\sigma_0,\sigma_1,\sigma_2\}$ be the Birman--Ko--Lee presentation
(see [2]) of the group of braids with three strings.
A {\it quasipositive factorization} of a braid $X\in\bold B_3$ is a collection
$(X_1,\dots,X_k)\in\bold B_3^k$ such that $X=X_1X_2\dots X_k$ and for each $i$, the braid $X_i$
is conjugate to $\sigma_1$. Note that $\sigma_0$, $\sigma_1$, and $\sigma_2$ are conjugate to
each other. We denote the set of quasipositive factorizations of $X$ by
$\Cal Q(X)$. A braid $X$ is called {\it quasipositive} if $\Cal Q(X)\ne\varnothing$.
The braid group $\bold B_k$ acts on $\Cal Q(X)$ by
$\Sigma_i:(X_1,\dots,X_k)\mapsto(Y_1,\dots,Y_k)$ where
$(Y_i,Y_{i+1})=(X_iX_{i+1}X_i^{-1},X_i)$ and $Y_j=X_j$ for $j\not\in\{i,i+1\}$.
This action is called the {\it Hurwitz action}. Elements belonging to the same orbit are
called {\it Hurwitz-equivalent}.

It is proven in [4] that each orbit of the Hurwitz action contains an element of
a certain explicitly specified finite set. The purpose of the present paper is to give
an easy criterion to decide if two given elements of this finite set belong to the same orbit.
To give precise statements, we need to introduce some notation with slightly differs from
that in [4].
Let us extend the alphabet $\Cal A$ up to $\hA=\Cal A\cup\{\hat\sigma_0,\hat\sigma_1,\hat\sigma_2\}$.
Let $\Cal A^*$ and
$\hA^*$ be the free monoids generated by $\Cal A$ and $\hA$ respectively.
If $U,V\in\hA^*$, then $U\letterwise V$ stands for equality in $\hA^*$
(i.~e.~letterwise coincidence of words) and $U=V$ (when $U,V\in\Cal A^*$) stands for
equality of the corresponding elements of $\bold B_3$.
For $U\in \hA^*$, we denote the word obtained from $U$ by erasing of all letters
$\hat\sigma_i$ (resp. by replacing each
$\hat\sigma_i$ with $\sigma_i$ or by replacing each $\sigma_i$ with $\hat\sigma_i$)
by $\bar U$ (resp. $U'$ or $\hat U$).
For example, if
$U\letterwise\sigma_0\hat\sigma_1\sigma_2\hat\sigma_1$, then
$U'\letterwise\sigma_0\sigma_2$,
$\bar U\letterwise\sigma_0\sigma_1\sigma_2\sigma_1$, and
$\hat U\letterwise\hat\sigma_0\hat\sigma_1\hat\sigma_2\hat\sigma_1$.

We set $\delta=\sigma_2\sigma_1=\sigma_1\sigma_0=\sigma_0\sigma_2$.
It is easy to check that any $X\in\bold B_3$ can be written as
$$
    X=U\delta^{-p},\qquad
    U\in\Cal A^*,\quad p\in\Bbb Z.                    \eqno(1)
$$
If, moreover, $U$ does not contain any subword which is equal (in $\bold B_3$) to $\delta$,
then the presentation of $X$ in the form (1) is unique and it is called
the {\it right Garside normal form}.

Let
$$
  W=W_1\hat x_1 W_2\hat x_2\dots W_k\hat x_k W_{k+1}\in \hA^*,
  \quad W_i\in\Cal A^*,
  \quad x_i\in\Cal A.                        \eqno(2)
$$
If $W'=\delta^p$, then $[W]$ stands for the quasipositive factorization of
$\bar W\delta^{-p}$ of the form
$(X_1,\dots,X_k)$ where $X_i=A_ix_iA_i^{-1}$, $A_i=W_1\dots W_i$.

To each $X\in\bold B_3$ we associate a graph $\Cal G_0(X)$ in the following way.
Let (1) be the right Garside normal form of $X$. We define the set of vertices of
$\Cal G_0(X)$ as
$\Cal V_0(X)=\{W\in\hA^*\mid W'=\delta^p, \bar W=U\}$.
Two vertices are connected by an edge if they are of the form
$A\hat x ByC$ and $AxB\hat yC$,
$x,y\in\Cal A$ where $xB'=B'y$ and, for some
$q\ge0$, either $B'=\delta^q$ (an edge of type (h1)), or
$xB'=\delta^q$ (an edge of type (h2)).

\proclaim{ Theorem 1 } Suppose that the right Garside normal form of a braid
$X\in\bold B_3$ is {\rm(1)} with $p\ge0$.
Then each orbit of the Hurwitz action on $\Cal Q(X)$ contains an element of the form
$[W]$, $W\in\Cal V_0(X)$. Two such elements belong to the same orbit if and only if
the corresponding vertices belong to the same connected component of $\Cal G_0(X)$.
\endproclaim

\noindent
{\bf Remark 1.} The condition $p\ge 0$ in Theorem 1 is not very restrictive. Indeed,
if $p<0$, then the Hurwitz action has a single orbit (see [4; Corollary 4]).

\medskip
The rest of the paper is devoted to the proof of Theorem 1.
Let $\tau:\hA^*\to\hA^*$ be the monoid homomorphism defined on the generators by
$\sigma_0\mapsto\sigma_1\mapsto\sigma_2\mapsto\sigma_0$,
$\hat\sigma_0\mapsto\hat\sigma_1\mapsto\hat\sigma_2\mapsto\hat\sigma_0$.
It induces the inner automorphism
$\tau:X\mapsto\delta^{-1}X\delta$ of $\bold B_3$.

\proclaim{ Lemma 1 } If $A,B\in\Cal A^*$, $A=B$, and $A\not\letterwise B$,
then $A=C\delta=\delta\tau(C)$, $C\in\Cal A^*$.
\endproclaim

\demo{ Proof }
It is known (see [2; Theorem 2.7]) that if $A=B$, then $B$ is obtained from
$A$ by applying the relations without inserting the inverses of the generators.
If $A\not\letterwise B$, then a relation can be applied to $A$, hence
$A=D\delta E$, $D,E\in\Cal A^*$.
Then $A=C\delta=\delta\tau(C)$ for $C=\tau(D)E$. \qed
\enddemo

\medskip\noindent
{\bf Definition.}
Let $W\letterwise a_1\dots a_n$, $a_i\in\hat{\Cal A}$, $W'=\delta^q$, $q>0$.
We say that letters $a_i$ and $a_j$ ($i<j$) match each other
in $W$ if $a_i,a_j\in\Cal A$,
$(a_i\dots a_j)'\letterwise a_iBa_j=\delta^{r+1}$
and $B=\delta^r$, $r\ge0$.
\medskip

Formally speaking, the indices $i$ and $j$ rather than the letters $a_i$ and $a_j$
match each other in the above definition. However, abusing the language, we shall
speak about matching letters to avoid cumbersome notation.

\proclaim{ Lemma 2 } {\rm(follows from [4; Lemma 5])}
Let $W\in\Cal A^*$, $W=\delta^q$, $q>0$. Then one can associate
parentheses to all letters of $W$ so that the parentheses are balanced
and matching pairs of parentheses correspond to matching pairs of letters.
\endproclaim

\proclaim{ Lemma 3 }
Let $W\letterwise AuBCvD\in\Cal A^*$, $W=\delta^q$, and
let $v$ match $u$ in $W$. Let $W_1\letterwise CvD\tau^q(AuB)$.
Then $W_1=\delta^q$ and $\tau^q(u)$ matches $v$ in $W_1$.
\endproclaim

\demo{ Proof } 
Let $E=AuB$, $F=CvD$. Then
$\delta^q=EF$, hence $W_1=F\delta^{-q}E\delta^q=F(F^{-1}E^{-1})E(EF)=EF=\delta^q$.
The same computation with $E=Au$, $F=BCvD$ yields
$BCvD\tau^q(Au)=\delta^q$. Since $BC=\delta^r$, we obtain
$vD\tau^q(Au)=\delta^{q-r}$. Similarly, by setting $E=A$ and $F=uBCvD$,
we obtain $D\tau^q(A)=\delta^{q-r-1}$.
\qed\enddemo

\proclaim{ Lemma 4 }
Let $W\letterwise ABC\in\Cal A^*$, $W=\delta^q$,
$B=\delta^r$. Then, for any letter $u$ in $A$ there is a letter in $A$ or in
$C$ which matches $u$ in $W$.
\endproclaim

\demo{ Proof } Follows from Lemma 3\ combined with
Lemma 2\ applied to $\tau^q(C)A$.
\qed\enddemo

Let a braid $X$ satisfies the hypothesis of Theorem 1.
We extend the graph $\Cal G_0(X)$ up to $\Cal G(X)$ as follows.
We define the set of vertices as
$\Cal V(X)=\bigcup_{j\ge0}\Cal V_j$ where $\Cal V_j=\{W\in\hA^*\mid
W'=\delta^{p+j}, \bar W\delta^{-p-j}=X\}$. The weight of a vertex $W\in\Cal V_j$
is defined as $\wt(W)=j$. 
We define edges of types (h1) and (h2) in the same way as in $\Cal G_0(X)$ and
we add more edges $(W,V)$ in the following cases:
\vbox{
\roster
\item"(h3)" $W\letterwise A\hat B_1 C$, $V\letterwise A\hat B_2 C$, 
$B_1,B_2\in \Cal A^*$, $B_1=B_2$;
\item"(v1)" $W\letterwise APB$, $V\letterwise A\tau^{-1}(B)$, 
            $P\in\Cal A^*$, $P=\delta$;
\item"(v2)" either
$W\letterwise APByC$ and $V\letterwise A\tau^{-1}(B\hat yC)$,
or $W\letterwise AyBPC$ and $V\letterwise A\hat yB\tau^{-1}(C)$
where $P\in\{u\hat x,\hat xu\}$, $\bar P=\delta$, and
$y$ matches $u$ in $W$;
\item"(v3)" $\wt(V)=\wt(W)-1$, $W$ is connected by an edge of type (h3) to
a vertex which is connected to $V$ by an edge of type (v2).
\endroster}
We call edges of types (h1--h3) {\it horizontal} and edges of types
(v1--v3) {\it vertical}.
The notation $(W,V)$ for vertical edges assumes that $\wt(W)>\wt(V)$.
If two vertices $W$ and $V$ belong to the same connected component of
$\Cal G(W)$, then we write $W\sim V$.
If either $W\letterwise V$ or there exists an edge $(W,V)$ of type, say, (h1), then
we write $W\sim_{\text{h1}}V$.
If $\wt W=\wt V$ and $W$ is connected to $V$ by a path in
$\Cal G(X)$, which passes through vertices of weight $\le\wt(W)$ only,
then we write $W\simh V$.

We define $f:(\hA\cup\Cal A^{-1})^*\to\hA^*$ by
$f(1)=1$, $f(xA)\letterwise xf(A)$ for $x\in\hA$ and
$f(x^{-1}A)\letterwise\tau(\tau(x)f(A))$ for $x\in\Cal A$.
Note that the braid $A^{-1}f(A)$ is a power of $\delta$.

\proclaim{ Lemma 5 }
(a). Let $W\in\Cal V(X)$. Then $W\sim f(Y_1\dots Y_k)$,
$Y_i\letterwise A_i\hat x_i A_i^{-1}$ for some
$A_i\in\Cal A^*$, $x_i\in\Cal A$,
$i=1,\dots,k$.
(b). If $X$ is conjugate to $\sigma_1$, then $\Cal V_0(X)$ contains a single element
and it has the form $f(AxA^{-1})$, $A\in\Cal A^*$, $x\in\Cal A$. 
\endproclaim

\demo{ Proof } (a).
Let $W$ be as in (2) and set $A_i\letterwise W_1\dots W_i$.
Then $f(Y_1\dots Y_k)$ is connected to $W$ by a chain of edges of type (v1).
(b). Follows from Lemma 2.
\qed\enddemo

\proclaim{ Lemma 6 } {\rm[4; Lemmas 6 and 7].}
If $W\in\Cal V(X)$ then $W\sim V$
for some $V\in\Cal V_0(X)$.
\qed
\endproclaim

\proclaim{ Lemma 7 } Let $W,V\in\Cal V(X)$. Then $[W]\sim[V]$
is equivalent to $W\sim V$.
\endproclaim

\demo{ Proof } $(\Rightarrow)$
Let $[W]=(X_1,\dots,X_i)$.
Lemma 5(a) combined with Lemmas 6\ and 5(b) applied to
$Y_i$'s imply
$W\sim f(Y_1\dots Y_k)$ where $\Cal V_0(X_i)=\{f(Y_i)\}$.
Thus $[W]=[V]\Rightarrow W\sim V$.
It remains to note that if $W\letterwise f(Y_1\dots Y_k)$,
$Y_i\letterwise A_i\hat x_i A_i^{-1}$, and
$$
V\letterwise
f(\dots Y_{i-1}\bar Y_i Y_{i+1}\bar Y_i^{-1}
 Y_i Y_{i+2}\dots ),\qquad
W_1\letterwise
f(\dots Y_{i-1} Y_i Y_{i+1}\bar Y_i^{-1}
 \bar Y_i Y_{i+2}\dots),
$$
then $\Sigma_i([W])=[V]$, $W_1\sim_{\text{h2}} V$, and $W_1$ is connected to $W$
by a chain of (v1)-edges.

\smallskip
$(\Leftarrow)$ (cp. [4; Lemma 6]).
It suffices to consider the cases
$W\sim_{\text{h1}}V$ and $W\sim_{\text{h2}}V$.
Let $W$ be as in (2). Then, for some $m$, $s$, $1\le m<s\le k+1$,
we have $W\letterwise A\hat x_m ByC$ and
$V\letterwise A x_m B\hat y C$ where $A\letterwise W_1\hat x_1W_2\hat x_2\dots W_m$,
$B\letterwise W_{m+1}\hat x_{m+1}\dots W_{s-1}\hat x_{s-1} D$,
$C\letterwise E\hat x_s W_{s+1}\hat x_{s+1}\dots \hat x_k W_{k+1}$
(if $s=k+1$, then $C\letterwise E$),
$W_s\letterwise DyE$.
Let $[W]=(X_1,\dots,X_k)$ and $[V]=(Y_1,\dots,Y_k)$.

It is clear that $Y_i=X_i$ for $i<m$.

If $m\le i < s-1$, then
$Y_i = B_i x_{i+1} B_i^{-1}$ where
$B_i=A'x_m W_{m+1}W_{m+2}\dots W_{i+1}$
$=(A'x_m(A')^{-1})A' W_{m+1}\dots W_{i+1}=X_m W_1\dots W_{i+1}$, whence
$Y_i = X_m X_{i+1} X_m^{-1}$.

If $i=s-1$, then $Y_i=A'x_m B' y (A'x_mB')^{-1}$. Using $x_m B'=B' y$, we obtain
$Y_i = A' x_m (x_m B')(A'x_mB')^{-1} = A' x_m (A')^{-1}=X_m$.

If $i\ge s$, then $Y_i = B_i x_i B_i^{-1}$ where
$B_i=A' x_m B' EF$ and
$F=W_{s+1}\dots W_i$. Using $x_m B'=B' y$,
we obtain $B_i = A'B'yEF = W_1W_2\dots W_i$, whence $Y_i=X_i$.

Thus, $[V] 
 = (X_1,\dots,X_{m-1},\,
    X_mX_{m+1}X_m^{-1},\dots, X_m X_{s-1} X_m^{-1},\, X_m,\,
     X_s,\dots,X_k)$
\par\noindent
$ = \Sigma_{s-2}\dots\Sigma_{m+1}\Sigma_m[W]$.
\qed
\enddemo

\break
Let $\bold c_s:\hA^*\to\hA^*$ be defined as
$aW\mapsto W\tau^s(a)$, $a\in\hA$, $W\in\hA^*$.

\proclaim{ Lemma 8 } {\rm(follows from Lemma 3)} (a).
Let $e=(W,V)\letterwise(aA,bB)$, $a,b\in\hA$,
be an edge of $\Cal G(X)$ of type (h1) or (h2). Let $s=p+\wt(W)$.
Then $e_1=(\bold c_s(W),\bold c_s(V))$ is an edge of $\Cal
G(\bold c_s(\bar W)\delta^{-s})$
of type (h1) or (h2). If $a=b$, then $e$ and $e_1$ are of the same type.
If $a\ne b$, then they are not.

\smallskip\noindent
(b). Let $e=(W,V)$ be an edge of $\Cal G(X)$ of type (v1) or (v2). Let $s=p+\wt(W)$.
Let $e_1=(\bold c^m_s(W),\bold c^m_{s-1}(V))$ where $m=2$ 
if $W$ starts with the word $P$ used in the definition of $e$, and
$m=1$ otherwise.
Then $e_1$ is an edge of $\Cal G(\bold c_s^m(\bar W)\delta^{-s})$ 
of the same type as $e$. \qed
\endproclaim

\proclaim{ Lemma 9\ \rm(Diamond Lemma)}
Let $e_1=(W,V_1)$ and $e_2=(W,V_2)$ be vertical edges in $\Cal G(X)$.
Then $V_1\simh V_2$.
\endproclaim

\demo{ Proof } 
By Lemma 8,  mutual arrangements of subwords of $W$
used in the definition of the edges can be considered up to cyclic
permutations.

Case 1. Both $e_1$ and $e_2$ are of type (v1). The statement is evident.

Case 2. Both $e_1$ and $e_2$ are of type (v2).
Let $P_i\in\{u_i\hat x_i,\hat x_iu_i\}$, and $y_i$, $i=1,2$, be the subwords of $W$
used in the definition of $e_i$.

Case 2.1. $P_1$ coincides with $P_2$, i.~e., $W\letterwise Ay_1By_2CP$,
$V_1\letterwise A\hat y_1By_2C$,
$V_2\letterwise Ay_1B\hat y_2C$, $P\in\{\hat xu,u\hat x\}$, $\bar P=\delta$.
By definition of edges of type (v2), we have $B'y_2C'=\delta^q$, $C'=\delta^r$,
$y_1B'y_2C'u=\delta^{q+1}$, $y_2C'u=\delta^{r+1}$. The former two identities
imply $B'y_2=\delta^{q-r}$. The latter two identities imply $y_1B'=\delta^{q-r}$.
Hence $V_1\sim_{\text{h2}}V_2$.

Case 2.2. $P_1$ and $P_2$ have a common letter $\hat x_1=\hat x_2=\hat x$.
Since $u_1x=xu_2=\delta$, it follows that $u_1$ cannot match $u_2$.
Hence $y_1$ and $y_2$ cannot coincide with them.

Case 2.2.1. $y_1=y_2=y$. $W\letterwise AyBu_1\hat x u_2$,
$V_1\letterwise A\hat y B\tau^{-1}(u_2)$, $V_2\letterwise A\hat yB u_1$.
Then $B'u_1=\delta^q$ and $B'=\delta^r$, whence $u_1=\delta^{q-r}$. Contradiction.

Case 2.2.2. $W\letterwise Ay_2By_1Cu_1\hat x u_2$,
$V_1\letterwise A y_2 B \hat y_1 C \tau^{-1}(u_2)$,
$V_2\letterwise A\hat y_2 B y_1 C u_1$.
Since $u_1\delta=u_1xu_2=\delta u_2$, we have $u_2=\tau(u_1)$,
whence $V_1\letterwise A y_2 B \hat y_1 C u_1$.
Since $B'y_1C'u_1=\delta^q$ ($e_2$ is of type (v2)) and $y_1C'u_1=\delta^{r+1}$
($e_1$ is of type (v2)), it follows that $B'=\delta^{q-r-1}$.
To conclude that $V_1\sim_{\text{h1}}V_2$, it remains to check that
$y_1=\tau^{q-r-1}(y_2)$. Indeed, $y_1\delta^r u_1=\delta^{r+1}$
implies $\tau^r(y_1)u_1=\delta$, and  $y_2\delta^q u_2=\delta^{q+1}$
implies $\tau^q(y_2)u_2=\delta$. Recall that $u_2=\tau(u_1)$. So, we obtain
$\delta=\tau(\delta)=\tau(\tau^r(y_1)u_1)
=\tau^{r+1}(y_1)u_2$,
whence $\tau^{r+1}(y_1)=\delta u_2^{-1}=\tau^q(y_2)$.

Case 2.2.3. $W\letterwise Ay_1By_2Cu_1\hat x u_2$,
$V_1\letterwise A\hat y_1 B y_2 C \tau^{-1}(u_2)$,
$V_2\letterwise A y_1 B\hat y_2 C u_1$.
Since $e_2$ is of type (v2), we have $C'u_1=\delta^q$.
Hence, by Lemma 2, we have $Cu_1\letterwise C_1z_1C_2u_1$
where $z_1$ matches $u_1$ in $Cu_1$.
Let $V_3\letterwise A y_1 B y_2 C_1 \hat z_1 C_2 u_1$.
Then $e_3=(W,V_3)$ is an edge of type (v2) and the pair $(e_1,e_3)$ (resp. $(e_3,e_2)$)
satisfies the conditions of Case 2.1 (resp. by Case 2.2.2).

Case 2.3. $P_1$ and $P_2$ have a common letter $u_1=u_2=u$,
i.~e., $W\letterwise A\hat x_1u\hat x_2$. Since $x_1\delta=x_1 u x_2=\delta x_2$,
we have $x_2=\tau^{-1}(x_1)$. Hence
the edge of type (v2) defined by $(P_1,y_2)$ coincides with $e_2$.
The problem is reduced to Case 2.1.

Case 2.4. $P_1$ and $P_2$ are disjoint, $y_1=u_2$, and $y_2=u_1$. Then
$W\letterwise A P_1 B P_2$, $V_1\letterwise A \tau^{-1}(B \hat P_2)$,
$V_2\letterwise A\hat P_1 B$, $B'=\delta^q$. If $q=0$, then
$V_1\sim_{\text{h3}}V_2$. Let $q>0$. Let
$B\letterwise\hat CzD$, $z\in\Cal A$ and
let $xy=yz=\delta$.
Since $P_1C=\tau^{-1}(C)xy$, we have
$V_3\letterwise A\tau^{-1}(\hat C)\hat x\hat y z D\sim_{\text{h3}}V_2$.
Then $V_1\letterwise E\tau^{-1}(zD \hat P_2)$ and
$V_3\letterwise E\hat x\hat y z D$ for $E\letterwise A\tau^{-1}(\hat C)$.
Since $zD'=(\hat CzD)'=B'=\delta^q$, it follows from Lemma 2\ that
$zD$ includes a letter
$u$ which matches $z$. Let $(V_3,U_3)$ 
be the edge of type (v2) defined by $(\hat y z,u)$.
Then $D\letterwise D_1uD_2$ and $U_3\letterwise E\hat x\tau^{-1}(D_1\hat u D_2)$. 
We have $D_1'=\delta^r$ and $zD'_1u=\delta^{r+1}$.
Recall that $zD_1'uD_2'=\delta^q$ which implies $D'_2=\delta^{q-r-1}$.
Let $v=\tau^{q-r-1}(u)$ and $vw=\delta$.
We may assume that $P_2\letterwise vw$ (otherwise we add an edge of type (h3)).
Then $V_1\sim_{\text{h1}}V_4\letterwise E\tau^{-1}(zD_1\hat u D_2 v\hat w)$.
Let us show that $\tau^{-1}(z)$ matches
$\tau^{-1}(v)$ in $V_4$. Indeed, we have
$D'_1D'_2=\delta^r\delta^{q-r-1}=\delta^{q-1}$ and $uD'_2=D'_2v$ (by the choice of $v$),
whence $zD'_1 D'_2v=zD'_1uD'_2=B'=\delta^q$.
Let $(V_4,U_4)$ be the edge of type (v2) defined by
$(\tau^{-1}(v\hat w),\tau^{-1}(z))$. Then
$U_4\letterwise E\tau^{-1}(\hat zD_1\hat uD_2)$. It remains to note that
$x\delta=xyz=\delta z$, whence $x=\tau^{-1}(z)$ and therefore $U_3\letterwise U_4$.

Case 2.5. $P_1$, $P_2$, and $y_2$ are pairwise disjoint:
$W\letterwise D_1 D_2$ where $D_1\letterwise A P_1 B$,
$D_2\letterwise y_2 C P_2$,
$C'=\delta^q$, and $D_2'=\delta^{q+1}$.
Then $D_1'=W'\delta^{-q-1}$ is a power of $\delta$. By
Lemma 2, this implies that $D_1$ has a letter $z_1$ which matches $u_1$.
Let $(W,V_3)$ be the edge of type (v2) defined by $(P_1,z_1)$.
Then we have $V_1\simh V_3$ (see Case 2.1),
$V_2\letterwise D_1 E_2$, and $V_3\letterwise E_1\tau^{-1}(D_2)$ where
$D_i\sim_{\text{v2}}E_i$ in the corresponding graphs.
Hence $V_2\sim_{\text{v2}} E_1\tau^{-1}(E_2)\sim_{\text{v2}} V_3$.

Case 3. $e_1$ is of type (v3) and $e_2$ is of type (v2) or (v3).
Let $W\sim_{\text{h3}}W_i\sim_{\text{v2}}V_i$, $i=1,2$.
Let $P_i\in\{\hat x_i u_i,u_i\hat x_i\}$ and $y_i$ be the subwords of $W_i$
used in the definition of the edge $(W_i,V_i)$.
We may assume that
$W_1\letterwise Eu_1\hat x_1\hat F_1 u_2$ and
$W_2\letterwise Eu_1\hat F_2\hat x_2 u_2$ where $x_1 F_1=F_2 x_2$
(otherwise the problem reduces to Case 2). By Lemma 1, we have
$x_1F_1=\delta\tau(D)=D\delta$, $D\in\Cal A^*$.
So, we may assume that
$x_1 F_1\letterwise x_1v_1\tau(D)$ and $F_2 x_2\letterwise D v_2 x_2$ where
$x_1 v_1=v_2 x_2=\delta$. Without loss of generality we may assume also that
$u_1x_1v_1\letterwise\sigma_2\sigma_1\sigma_0$.

Case 3.1. $y_1=u_2$ and $y_2=u_1$. Then $u_1(\hat x_1\hat F_1)'u_2=u_1u_2=\delta$.
Hence
$W_1\letterwise E \sigma_2\hat\sigma_1\hat\sigma_0\tau(\hat D)\sigma_1$ and
$W_2\letterwise E \sigma_2\hat D\hat\sigma_0\hat\sigma_2\sigma_1$, whence
$V_1\letterwise E\tau^{-1}(\hat\sigma_0\tau(\hat D)\hat\sigma_1)
\letterwise E\hat\sigma_2\hat D\hat\sigma_0\letterwise V_2$.

Case 3.2. $E\letterwise E_2y_1E_1$. Then $E'_1=\delta^q$ and
$y_1E'_1u_1=\delta^{q+1}$. Hence, by Lemma 4, there is a letter $z_2$
in $E_2$ which matches $u_2$. Let $(W_2,V_3)$ be the edge
of type (v2) defined by $(P_2,z_2)$. Then
$V_2\simh V_3$ (see Case 2.1).
By passing to a cyclic permutation of $W$, we may assume that
$W_1\letterwise Ay_1Bu_1\hat x_1\hat v_1\tau(\hat D)u_2Cz_2$,
$W_2\letterwise Ay_1Bu_1D\hat v_2\hat x_2u_2Cz_2$,
$B'=\delta^q$, $y_1Bu_1=\delta^{q+1}$, $C'=\delta^r$, $u_2C'z_2=\delta^{r+1}$.
For each possible value of $u_2$,
we find $V_4$ such that $V_1\simh V_4\simh V_3$.
(for $u_2=\sigma_1,\sigma_2$ we just list the indices which should be inserted into
the formulas for $V_1$, $V_3$, and $V_4$).
$$
\split
 u_2&=\sigma_0:
\\
 V_1&\letterwise A\tau^{-q}(\hat\sigma_0)B\tau^{-1}\big(\hat\sigma_0
   \tau(\hat D) \sigma_0C\tau^r(\sigma_2)\big)
   \letterwise
   A\tau^{-q}(\hat\sigma_0)B \hat\sigma_2
   \hat D \sigma_2\tau^{-1}(C)\tau^r(\sigma_1),
\\
 V_3&\letterwise A\,\tau^{-q}(\sigma_0)\,B\,\sigma_2 \,\hat D\, 
     \hat\sigma_2\,\tau^{-1}\big( C\,\tau^r(\hat\sigma_2)\big)
   \,\letterwise
    A\tau^{-q}(\sigma_0)B\sigma_2 \hat D 
     \hat\sigma_2\tau^{-1}(C)\tau^r(\hat\sigma_1),
\\
 V_4&\letterwise
   A\tau^{-q}(\sigma_0)B\hat\sigma_2 \hat D
         \sigma_2\tau^{-1}(C)\tau^r(\hat\sigma_1);
\\
u_2&=\sigma_1{:}\,\;
V_1(\hat0 \hat2 0 2),\; V_4(0 \hat2 \hat0 2),\; V_3(0 2 \hat0 \hat2);
 \quad
u_2=\sigma_2{:}\,\;
V_1(\hat0 \hat2 1 0),\; V_4(\hat0 2 1 \hat0),\; V_3(0 2 \hat1 \hat0).
\endsplit
$$

Case 4. $e_1$ is of type (v1) and $e_2$ is of type (v2).
Let $(W,V_1)\letterwise (AP_1,A)$, $P_1=\delta$.
Let $P$ and $y$ be the subwords of $W$ used in the definition of $e_2$.

Case 4.1. $P_1$ has one common letter with $P$:
$W\letterwise AyB\hat xuv$, $xu=uv=\delta$, $B'=\delta^q$, $yB'u=\delta^{q+1}$,
$V_1\letterwise AyB\hat x$, and $V_2\letterwise A\hat yB\tau^{-1}(v)$.
Then $x\delta=xuv=\delta v$ implies $\tau^{-1}(v)=x$, and
$y\delta^qu=\delta^{q+1}$ combined with $xu=\delta$ implies $x=\tau^q(y)$.
Hence $V_1\sim_{\text{v1}}V_2$.

Case 4.2. $P_1$ is disjoint from $P$.
Lemma 4\ implies $A\letterwise BzCPD$ where $z$ matches $u$.
Let $(W,V_3)$ be the edge of type (v2) defined by $(P,z)$.
Then $V_2\simh V_3$ (see Case 2.1)
and there exist vertical edges $(V_1,V_4)$ and $(V_3,V_4)$
where $V_4\letterwise B\hat zC\tau^{-1}(D)$.

Case 5. $e_1$ is of type (v1) and $e_2$ is of type (v3).
Let $W\sim_{\text{h3}}W_1\sim_{\text{v2}}V_2$.
Then we have $W_1\sim_{\text{v1}}V_1$ which reduces the problem to Case 4.
\qed
\enddemo

\proclaim{ Lemma 10 }
Let $e=(W_1,W_2)$ be a horizontal edge in $\Cal G(X)$. If
$\wt(W_1)>0$, then there exist vertical edges $(W_1,V_1)$ and $(W_2,V_2)$
such that $V_1\simh V_2$.
\endproclaim

\demo{ Proof } The condition $\wt(W_1)>0$ implies that
$\bar W_1$ has a subword $P$ which is equal to $\delta$.
By Lemma 8, it suffices to consider edges of type (h1) and (h3) only.

Case 1. $e$ of type (h1). Let $x$ and $y$ be as in the definition of edges of type (h1).
It is clear that $P\not\letterwise xy$. If $P$, $x$, $y$ are disjoint,
then the statement is evident.

Case 1.1. $W_1\letterwise A\hat y B ux$,
$W_2\letterwise AyBu\hat x$, $ux=\delta$. 
Let $V_1=A\hat yB$. Then $W_1\sim_{\text{v1}}V_1$. Since $e$ is an edge
of type (h1), we have $B'u=\delta^q$, $x=\tau^q(y)$. Hence, by Lemma 2,
there is a letter $z$ in $Bu$ which matches $u$, i.~e., $B\letterwise B_1zB_2$,
$B'_2=\delta^r$, $zB'_2u=\delta^{r+1}$. Hence
$W_2\sim_{\text{v2}} V_2\letterwise AyB_1\hat z B_2$.
Let us show that $V_1\sim_{\text{h1}}V_2$. Indeed, $B'_1zB'_2u=B'u=\delta^q$
combined with $zB_2'u=\delta^{r+1}$ imply $B_1'=\delta^{q-r-1}$, 
and $z\delta^r u=\delta^{r+1}$ combined with $ux=\delta$ imply $x=\tau^{r+1}(z)$.
Since $x=\tau^q(y)$, we obtain $z=\tau^{q-r-1}(y)$.

Case 1.2. $W_1\letterwise A\hat y B xu$, $W_2\letterwise AyB\hat xu$, $xu=\delta$,
$B'=\delta^q$, and $yB'=B'x$. Then 
$yB'u=B'(xu)=\delta^{q+1}$, i.~e., $(W_1,A\hat yB)$ and $(W_2,A\hat yB)$
are edges of type (v1) and (v2) respectively.

Case 2. $e$ of type (h3). Without loss of generality we may assume that
$W_1\letterwise Au\hat B_1$, $W_2\letterwise Au\hat B_2$,
$u\in\Cal A$, $B_1,B_2\in\Cal A^*$, and $B_1=B_2$.
Lemma 1\ implies $B_1=\delta B$. Let $W\letterwise Au\hat x\hat v\hat B$
where $ux=xv=\delta$. 
Then $W_1\sim_{\text{h3}}W\sim_{\text{h3}}W_2$.
Let $y$ matches $u$ and let $(W,V)$ be the edge of type (v2)
defined by $(u\hat x,y)$.
Then $W_1\sim_{\text{v3}}V\sim_{\text{v3}}W_2$.
\qed
\enddemo

\demo{ Proof of Theorem 1}
By Lemma 7, it is enough to show that:
if $W\sim V$ and $\wt(W)=\wt(V)$, then
$W\simh V$. Indeed,  if a path from $W$ to $V$
in $\Cal G(X)$ passes through vertices whose weight is greater than $\wt(W)$, then,
by Lemmas 9\ and 10, we can modify the path so that either
the number of vertices of the maximal weight is reduced, or
the number of horizontal edges incident to such vertices is reduced.
\enddemo

\smallskip
ACKNOWLEDGMENTS.
This work was supported by the Russian Science
Foundation, project no. 14-21-00053.


\Refs
\def\r{\ref}

\r\no1
\by     E.~Artal Bartolo, J.~Carmona Ruber, J.~I.~Cogolludo Agust{\'\i}n
\paper  Effective invariants of braid monodromy
\jour   Trans. Amer. Math. Soc. \vol 359\yr 2007 \pages 165--183
\endref

\r\no2
\by     J.~Birman, K.-H.~Ko, S.-J.~Lee
\paper  A new approach to the word and conjugacy problems in the braid groups
\jour   Adv. Math. \vol 139 \yr 1998 \pages 322--353
\endref

\r\no3
\by     Vik.~S.~Kulikov
\paper  Factorization semigroups and irreducible components of the Hurwitz space. II
\jour   Izv. Math., \vol 76:2 \yr 2012 \pages 356--364 \endref

\r\no4
\by     S.~Yu.~Orevkov
\paper  On the Hurwitz action on quasipositive factorizations of 3-braids
\jour   Doklady Math. \vol 91 \yr 2015 \pages 173--177
\endref

\endRefs
\enddocument